\documentclass[a4paper,11pt]{article}
\usepackage{indentfirst,latexsym,bm,amsmath,amssymb,amsthm}
\usepackage[dvips]{graphicx}

\makeatletter\@addtoreset{equation}{section} \makeatother
\newtheorem{theorem}{Theorem}[section]
\newtheorem{lemma}[theorem]{Lemma}

\newtheorem{remark}{Remark}[section]

 


\begin{document}
\date{}
\title{ Blow up for some semilinear wave equations in multi-space
dimensions }

\author{Yi Zhou\thanks{  School of Mathematical Sciences, Fudan
University, Shanghai 200433, P. R. China ({\tt Email:
yizhou@fudan.ac.cn})
  } \and
Wei Han\thanks{ Department of Mathematics, North University of
China,  Taiyuan, Shanxi  030051, P. R. China  ({\tt Email:
sh\_hanweiwei1@126.com})}.}
\date{}

\maketitle
\begin{abstract}

 In this paper, we discuss a new nonlinear phenomenon.  We find that in $n\geq 2$ space dimensions, there exists two indexes $p$ and $q$ such that the cauchy problems for the nonlinear wave equations
  \begin{equation} \label{0.1}
        \Box u(t,x ) = |u(t,x)|^{ q}  , \  \    x\in R^{n},
        \end{equation}
 and
     \begin{equation} \label{0.2}
      \Box u(t,x ) =    |u_{t}(t,x)|^{ p} , \  \    x\in R^{n}
     \end{equation}
 both have global existence for small initial data, while  for the combined
  nonlinearity,   the solutions to the Cauchy problem for the nonlinear wave equation
   \begin{equation} \label{0.3}     \Box u(t,x ) =  |  u_{t}(t,x)|^{ p} + |u(t,x)|^{ q}, \  \    x\in R^{n},
   \end{equation}
 with small initial data will blow up in finite time.   In the two dimensional case , we also find that if $  q=4$,  the Cauchy problem for the equation \eqref{0.1} has global existence,  and the Cauchy problem for the equation
   \begin{equation} \label{0.4}
      \Box u(t,x ) =    u (t,x)u_{t}(t,x)^{ 2} , \  \    x\in R^{2}
     \end{equation}
  has almost global existence, that is, the life span is at least $ \exp (c\varepsilon^{-2}) $ for initial data of size $ \varepsilon$. However, in the combined
  nonlinearity case,   the Cauchy problem for the equation
  \begin{equation} \label{0.5}     \Box u(t,x ) =  u(t,x)  u_{t}(t,x)^{ 2} + u(t,x)^{ 4}, \  \    x\in R^{2}
   \end{equation}
  has a life span which is  of the order of  $  \varepsilon^{ -18} $ for the initial data of size $ \varepsilon$,  this is considerably shorter in magnitude than that of the first two equations.
  This solves an open optimality
 problem for general theory  of fully nonlinear wave equations  (see \cite{Katayama}).

\par {\bf Keywords:} Fully nonlinear wave equations ;  Life-Span;
Cauchy problem

\end{abstract}
\section{Introduction and Main Results}

  First we shall outline the general theory on the Cauchy problem for
 the following $n$-dimensional fully  nonlinear wave
 equations:
\begin{equation} \label{1.1}
\left \{
\begin{array}{lllll}
u_{tt}- \Delta u= F(u, Du, D_{x}Du ),  \   \    (\Delta
=\frac{\partial^{2} }{\partial x_{1}^{2}}+ \cdots \frac{\partial^{2}
}{\partial x_{n}^{2}}   ),
 \         \cr\noalign{\vskip2mm}
 t=0:   \    \    \       u=\varepsilon f(x),  \   \          u_{t}=\varepsilon g(x),   \    \
   ( x=(x_{1}, \cdots, x_{n})) ,
\end{array} \right.
\end{equation}
where $$ Du= (u_{x_{0}},
u_{x_{1}}, \cdots, u_{x_{n}}),  \    \       x_{0}=t, $$
   $$ D_{x}Du=( u_{x_{i}x_{j}}, \    i,j=0,1,\cdots,n, i+j\geq 1),   $$
   $ f(x), g(x) \in C^{\infty}_{ 0} (R^{n} ) $ and $\varepsilon>0 $
   is a small parameter. Here, for simplicity of notations we write
   $ x_{0}=t $.

Let
$$ \hat{\lambda} =\big( \lambda; \   (\lambda_{i}),i=0,1,\cdots,n; \   (\lambda_{ij}), i,j=0,1,\cdots, n,
 \  i+j\geq 1  \big).  $$
Suppose that in a neighborhood of $ \hat{\lambda}=0 $, say, for $
|\hat{\lambda}|\leq 1, $  the nonlinear term $F=F(\hat{\lambda}) $
in equation  \eqref{1.1} is a sufficiently smooth function with
$$ F(\hat{\lambda})= O(|\hat{\lambda}|^{ 1+ \alpha} ), $$
 where $ \alpha$ is an integer and $\alpha \geq 1. $

We define a lifespan $T(\varepsilon)$ of  solutions to problem
\eqref{1.1}  to be the supremum of all  $\tau>0 $ such that there
exists a classical solution to \eqref{1.1} for $x\in R^{2}$ on
$0\leq t< T(\varepsilon)$.  When  $ T(\varepsilon)=+ \infty$, we mean that the problem \eqref{1.1} has global existence.

In chapter 2 of Li and Chen  \cite{Ta-tsien}, we have long histories
on the estimate for  $ T(\varepsilon) $.   The lower bounds of  $
T(\varepsilon) $  are summarized in the following table.  Let  $
a=a( \varepsilon ) $  satisfy
  $$   a^{ 2} \varepsilon^{ 2} \log (a+1) =1   $$
 and $c$ stand for a positive constant independent of $ \varepsilon
 $.  We have(see also a table in Li  \cite{Ta-tsien 4})
\begin{table}[ht]
\caption{ General Theory for the sharp lower bound of the lifespan
for fully nonlinear wave equations
 } 
\centering 

\begin{tabular}{|c| c |c |c|} 
\hline\hline 
 $ T( \varepsilon ) \geq $ &  $  \alpha=1  $&   $ \alpha=2 $ &  $  \alpha\geq 3  $ \\ [0.5ex] 
\hline 
 n=2 &    $ \begin{array}{lllll}
 c a( \varepsilon ),  \   \   \   \   \mbox{ in general case,}
 \       \cr\noalign{\vskip3mm}
  c\varepsilon^{ -1 },    \   \     \   \    \mbox{ if }  \displaystyle\int_{ R^{ 2}} g(x)
  dx=0,  \cr\noalign{\vskip3mm}
 c\varepsilon^{ -2 },    \   \     \   \     \mbox{ if } \partial_{ u}^{ 2}
 F(0) =0
\end{array}  $      &
 $ \begin{array}{lllll}
  c\varepsilon^{ -6 },   \   \   \   \   \mbox{ in general case,}
 \       \cr\noalign{\vskip3mm}
  c\varepsilon^{ -18 },    \   \     \   \    \mbox{ if }   \partial_{ u}^{ 3}
 F(0) =0
 \cr\noalign{\vskip3mm}
 \exp ( c\varepsilon^{ -2 }),    \   \     \   \    \mbox{ if }   \partial_{ u}^{ l}
 F(0) =0 ( l=3, 4 )
\end{array}  $
 &  $ \infty $ \\ 

\hline

 n=3 &    $ \begin{array}{lllll}
  c\varepsilon^{ -2 },   \   \   \   \   \mbox{ in general case,}
 \       \cr\noalign{\vskip4mm}
 \exp ( c\varepsilon^{ -1 }),    \   \     \   \    \mbox{ if }   \partial_{ u}^{ 2}
 F(0) =0
\end{array}  $
 & $ \infty $ & $ \infty $ \\
 \hline
n=4 &  $ \begin{array}{lllll}
  \exp ( c\varepsilon^{ -2 }),   \   \   \   \   \mbox{ in general case,}
 \       \cr\noalign{\vskip4mm}
  \infty,    \   \     \   \    \mbox{ if }   \partial_{ u}^{ 2}
 F(0) =0
\end{array}  $ &  $ \infty $ & $ \infty $  \\
\hline
 $ n\geq 5 $  &  $ \infty $  & $ \infty $ & $ \infty $ \\
 [1ex] 
\hline 
\end{tabular}

\label{table:nonlin} 
\end{table}

 We note that  all these lower bounds are known  to be sharp except for the case $ (n, \alpha)    = ( 2, 2)  $
 and $\partial_{ u}^{ 3}F(0) =0$. The aim of this paper is to show that in this case the lower bound obtained by \cite{Katayama} is indeed sharp. Therefore, we solve an open optimality
 problem for general theory  of fully nonlinear wave equations. We remark that the sharpness for
$  (n, \alpha)    = ( 4, 1)$  was only recently proved by
 Takamura and Wakasa \cite{Takamura and Wakasa}, see also Zhou and Han
\cite{Y. Zhou 5}.  For the  case $ (n, \alpha)    = ( 2, 2)  $ and  $\partial_{ u}^{ l} F(0) =0 ( l=3, 4 )$, the sharpness is due to Zhou and Han \cite{Y. Zhou 4},  in this case, it was believed  that  $ l=4 $ is a technical condition which may be
   removed,  however, we show in this paper that it is not the case,
    we show that if we drop this condition, the lifespan will be
    much shorter.   Therefore our result shows that the condition $ l=4 $ is necessary.

In this paper, we firstly consider the following Cauchy problem with
small initial data in two space dimensions
\begin{equation} \label{1.2}
\left \{
\begin{array}{lllll}
  \Box u(t,x ) =  u(t,x)u_{t} (t,x)^{ 2} + u^{ 4} (t,x),
  \   \   \
  x\in R^{ 2},  \   \    t>0,           \\
  t=0:  \     u=0,  \   \    \     \      u_{t}=\varepsilon g(x),   \  \
   \   x\in R^{ 2},
\end{array} \right.
\end{equation}
where   $ \Box= \partial_{t}^{ 2} -\sum\limits_{ i=1}^{ n}
\partial_{ x_{i}}^{ 2}$ is the wave operator,  and  $  g\in C^{ \infty}_{
0} ( R^{ n} )$,      $ \varepsilon >0 $ is a small parameter.

   For problem   \eqref{1.2},       what
interesting about this problem is that the Cauchy problem for the equation
  $$      \Box u(t,x ) = u^{ 4} (t,x), \   \   \
  x\in R^{ 2},  \   \    t>0            $$
  has global existence(see \cite{Glassey 2}),
     and the Cauchy problem for the equation
     $$   \Box u(t,x ) =    u(t,x)u_{t}^{ 2} (t,x), \   \   \
  x\in R^{ 2},  \   \    t>0              $$
  has almost global existence(see \cite{Ta-tsien}).   However,  in the combined
  nonlinearity case, the Cauchy problem for the equation
  $$    \Box u =  u  u_{t}^{ 2} + u^{ 4}, \   \   \
  x\in R^{ 2},  \   \    t>0   $$
   has a life span which is  of the order  $  \varepsilon^{ -18} $,  this is considerably shorter in magnitude  than that of
     $$     \Box u =    u(t,x)u_{t}(t,x)^{ 2}   $$
     and
       $$     \Box u =    u^{ 4} (t,x ) .   $$

For problem \eqref{1.2},  we consider compactly supported, radial,
nonnegative data  $ g\in C^{\infty}_{0}(R^{2})$,  and satisfy
\begin{equation}\label{1.4}
   \    g(x)=g(|x|)> 0, \  \mbox{ for}  \   |x|<1      \  \   g(x)=0, \
 \   \mbox{ for}  \   |x|>1.   \  \
\end{equation}

  We establish the following theorem for \eqref{1.2}:

\noindent\begin{theorem}\label{thm:1.1}   Let  $g $ is a smooth
function with compact support $g\in C_{ 0}^{ \infty} (R^{2}) $ and
satisfies \eqref{1.4},
         space dimensions $n=2 $. Suppose
  that problem   \eqref{1.2} has a solution $(u, u_{t})\in
C([0,T), \ H^{1}(R^{2})\times L^{3}(R^{2}))$ such that
 $$\emph{supp}(u, u_{t})\subset \{ (x,t): \    |x|\leq t+1    \}.$$
      Then   the solution $ u=u(t,x) $  will blow up in finite time, that is  $T<\infty$.  Moreover, we have the following estimates for the
   lifespan $ T(  \varepsilon) $ of solutions of  \eqref{1.2}:
        there exists a
 positive constant $A $ which is independent of $ \varepsilon$
 such that
\begin{equation}\label{1.5}
\begin{array}{ll}
  T(\varepsilon) \leq
  A\varepsilon^{-18  }.
\end{array}
\end{equation}
\end{theorem}

Secondly,  we will consider the following Cauchy problem with small initial data in $n( n\geq 2)$ space dimensions
\begin{equation} \label{1.3}
\left \{
\begin{array}{lllll}
  \Box u(t,x ) =  |  u_{t} (t,x)|^{ p} +  |u(t,x) |^{ q},
  \   \   \
  x\in R^{ n},  \   \    t>0,           \\
  t=0:  \     u=\varepsilon f(x),  \   \    \     \      u_{t}=\varepsilon g(x),   \  \
   \   x\in R^{ n},
\end{array} \right.
\end{equation}
where $ f, g\in C^{ \infty}_{ 0} ( R^{ n} )$,      $ \varepsilon >0 $ is a small parameter,   $
p>1, q>1  $.

For problem \eqref{1.3},  we consider compactly supported
nonnegative data  $f, g\in C^{\infty}_{0}(R^{n})$, $n \geq 2$ and
satisfy
\begin{equation}\label{1.6}
 f(x)\geq 0 , \    g(x)\geq 0, \      \  \   \ f(x)=g(x)=0, \
 \   \mbox{ for}  \   |x|>1 \   \
 \          \mbox{and} \  \    g(x)  \not\equiv
 0.
\end{equation}

  We establish the following theorem for \eqref{1.3}:

\noindent\begin{theorem}\label{thm:1.2}   Let  $f,g $ are smooth
functions with compact support $f,g\in C_{ 0}^{ \infty} (R^{n }) $
and satisfy \eqref{1.6},
         space dimensions $n \geq 2 $. Suppose
  that problem  \eqref{1.3}  has a solution $(u, u_{t})\in
C([0,T), \ H^{1}(R^{n})\times L^{r}(R^{n}))$,  where  $ r=\max(2, p) $ such that
 $$\emph{supp}(u, u_{t})\subset \{ (x,t): \    |x|\leq t+1    \},
 $$
 and the index $p, \   q $  satisfies the following conditions:

  \begin{equation}\label{1.7}
        \max\left(1, \frac{2 }{n-1}\right)<p\leq  \frac{  2n}{ n-1  };
    \end{equation}

\begin{equation}\label{1.8}   1<q<  \min\left(  \frac{  4 }{ ( n-1 )p -2  } +1, \    \frac{2n}{n-2} \right)     .
   \end{equation}
     Then   the solution $ u=u(t,x) $  will blow up in finite time, that is  $T<\infty$.  Moreover, we have the following estimates for the
   lifespan $ T(  \varepsilon) $ of solutions of   \eqref{1.3} :
        there exists a
 positive constant $A $ which is independent of $ \varepsilon$
 such that
\begin{equation}\label{1.9}
\begin{array}{ll}
  T(\varepsilon) \leq
  A\varepsilon^{-  \frac{ 2p(q-1)   }{ 2q+2- (n-1) p(q-1)  }  }.
\end{array}
\end{equation}
\end{theorem}

\begin{remark}\label{Remark:1.1}
In the Theorem \ref{thm:1.2}, we restrict $ q< 2^{*}=\frac{ 2n}{n-2}  $  in  the condition \eqref{1.8}, just to make that $ H^{1}(R^{n})\hookrightarrow L^{q}(R^{n})$ , so the nonlinearity  $ |u|^{q} $ can be integrable in $R^{n} $.
\end{remark}

\begin{remark}\label{Remark:1.2}
If we take $n=2,$ $p=3$ and $q=4$ in theorem 1.2, we can also obtain the upper bound of the life span  $A\varepsilon^{-18}$.
\end{remark}

\begin{remark}\label{Remark:1.3}
   In problem  \eqref{1.3},
   let $ p_{0}= p_{0}(n)= \frac{ 2}{ n-1} +1  $, and  $ q_{0 } =q_{
0}(n ) $  is the positive root of the quadratic equation $ \gamma(n, q)= (n-1) q^{
2} -( n+1 )q-2=0 $, that is
 $$  q_{ 0} =  q_{ 0} (n) =  \frac{  n+1+  \sqrt{ n^{ 2} +10n-7     }   }{  2( n-1 )   } .    $$
 Then we know that  $ p_{ 0} $ and $ q_{ 0} $ are the critical index of the following
semilinear wave equations respectively:

  \begin{equation} \label{4.25}      \Box u(t,x ) =  |  u_{t} (t,x)|^{ p} ,
    \end{equation}
and
  \begin{equation} \label{4.26}    \Box u(t,x ) =   |u(t,x) |^{ q}.
    \end{equation}
If $ p>p_{0} $, then the solution of the initial problem for the
above equation \eqref{4.25} will exists globally,
see \cite{Hidano, Y. Zhou 2};
    and also if $
q>q_{0} $, then the solution of the initial problem for the above
equation \eqref{4.26} will exists globally, this problem has a long history, one can see \cite{John 1, Glassey 1, Glassey 2, Georgiev, Sideris, Y. Zhou 1}.
 But  it can be showed that there exisits $ p>p_{0} $ and  $ q>q_{0} $  such that the conditions
 \eqref{1.7} and \eqref{1.8} can be still satisfied, for that purpose, we take for example $p>p_0$ and sufficiently close to $p_0$, then \eqref{1.7} will be satisfied. Take $q$ satisfies \eqref{1.8} and sufficiently close to $ \frac{  4 }{ ( n-1 )p -2 } +1  $ then $q$ will be sufficiently close to $ \frac{  4 }{  n-1  } +1  $ which is larger than $q_0$, this can be seen from the fact that
 \begin{equation} \label{4.27}
    \gamma\left(n, \frac{4}{n-1}+1   \right)=\frac{ 8}{n-1}>0.
  \end{equation}
Consequently,  from Theorem  \ref{thm:1.2},   the solutions of the Cauchy problem \eqref{1.3}
   will blow up in finite time, while the Cauchy problems for \eqref{4.25} and \eqref{4.26} have global existence.
\end{remark}

The rest of the paper is arranged as follows. We state  a
 preliminary Lemma in Section 2.   In Section 3, we prove the sharpness of the lower bound obtained by \cite{Katayama},  i.e., Theorem 1.1.
  Section 4 is devoted to the proof for our Theorem \ref{thm:1.2},  i.e., blow up and the upper bound estimate of lifespan
  of solutions  to some semilinear wave equations with small initial data in $n(n\geq 2)$  space dimensions.

\section{Preliminaries} \vskip .5cm

To prove the main results in this paper, we will employ the
following important ODE result:

\noindent \begin{lemma}\label{thm:2.1}  (see \cite{Sideris}, also see \cite{Y. Zhou 3}) Let
$ \beta>1$, $a\geq 1$, and
 $(\beta-1)a> \alpha-2$.  If $F\in C^{2}([0, T))$ satisfies
\\ (1) \   $F(t)\geq \delta (t+1)^{a}$,
 \\  (2) \   $\frac{d^{2}F(t)}{dt^{2}}\geq k (t+1)^{- \alpha}[F(t)]^{\beta }$,
\\     with some positive constants $\delta $, $k$,
then $F(t)$ will blow up in finite time,  $T<\infty$.  Furthermore,
we have the the following estimate for the life span $T(\delta)$ of
$F(t)$ :
   \begin{equation} \label{2.1}
   T(\delta)\leq c \delta^{ -\frac{(\beta-1)}{ (\beta-1)a-\alpha+2}
   },
    \end{equation}
     where  $c$ is a positive constant depending on $k$  but independent of $\delta$.
\end{lemma}

\noindent\begin{proof}\   For the proof of blow up result part see Sideris \cite{Sideris}. For the estimate of the life span of $F(t)$,  one can see Lemma 2.1 in \cite{Y. Zhou 3}.
\end{proof}

To outline the method, we following  Yordanov and Zhang  \cite{Q. S.
Zhang 1} will introduce the following functions: Let
     $$      \phi_{1}(x) =
     \displaystyle\int_{S^{n-1}} e^{x\cdot \omega} d\omega \geq 0  . $$
Obviously,  $ \phi_{1}(x)$  satisfy:  $  \Delta \phi_{1}= \phi_{1}
 $.

When space dimensions $ n\geq 2 $,    by rotation invariance, we
have
$$\begin{array}{ll}
     \phi_{1}(x) & =   \displaystyle\int_{S^{n-1}} e^{ | x|  \cdot \omega_{ 1}} d\omega
    \cr\noalign{\vskip 4mm} &
      =  \displaystyle\int_{ \omega_{ 1}^{ 2} + \widetilde{\omega}^{ 2} =1   } e^{ | x|  \cdot \omega_{ 1}} d\omega_{
      1} d\widetilde{\omega}  \cr\noalign{\vskip 4mm} &
      =  c \displaystyle\int_{0 }^{1} e^{ | x|  \cdot \omega_{ 1}} ( \sqrt{ 1- \omega_{ 1}^{ 2} } )^{
       n-3} d \omega_{ 1}  + c \displaystyle\int_{0 }^{1} e^{ -| x|  \cdot \omega_{ 1}} ( \sqrt{ 1- \omega_{ 1}^{ 2} } )^{
       n-3} d \omega_{ 1}\cr\noalign{\vskip 4mm} &
      \le  c e^{ | x| } \displaystyle\int_{0 }^{1} e^{ -| x|  \cdot ( 1- \omega_{ 1} )
      } ( \sqrt{ 1- \omega_{ 1}^{ 2} } )^{
       n-3} d \omega_{ 1}  +C\cr\noalign{\vskip 4mm} &
       \leq   c e^{ | x| } \displaystyle\int_{0 }^{1} e^{ -| x|  \cdot ( 1- \omega_{ 1} )
      } (  1- \omega_{ 1} )^{
       \frac{n-3}{ 2}} d \omega_{ 1} +C  \   \      ( \mbox{ we  take }   \lambda =| x| ( 1- \omega_{ 1} ) )
       \cr\noalign{\vskip 4mm} &
  =  c e^{ | x| } | x|^{ -\frac{ n-1}{ 2 } } \displaystyle\int_{0 }^{ | x|
  } e^{- \lambda } \lambda^{ \frac{ n-3}{ 2 } } d\lambda  +C\cr\noalign{\vskip 4mm} &
   \leq   c e^{ | x| } | x|^{ -\frac{ n-1}{ 2 } } \displaystyle\int_{0 }^{  \infty
  } e^{- \lambda } \lambda^{ \frac{ n-3}{ 2 } } d\lambda  +C\cr\noalign{\vskip 4mm} &
   = \widetilde{C} e^{ | x| } | x|^{ -\frac{ n-1}{ 2 } }   .
\end{array}$$

Moreover,  obviously we have

 $$   | \phi_{1}(x) | \leq   e^{ |x| }  \displaystyle\int_{S^{n-1}} d\omega = C e^{ |x| }.                $$
  Thus we can conclude that
  \begin{equation} \label{2.2}
   | \phi_{1}(x) | \leq  C e^{ |x| } ( 1+ | x|)^{ -\frac{ n-1}{ 2 } }   \    \       (  n\geq  2).
    \end{equation}

By the positivity of $\phi_{1}(x)  $,  so we get that when $ n\geq 2
$,

$$   0\leq   \phi_{1}(x) \leq  C e^{ |x| } ( 1+ | x|)^{ -\frac{ n-1}{ 2 } } ,   \   \      (  n\geq  2)  .            $$

In order to describe the following  methods, we define the following
test function
\begin{equation} \label{2.4}
\psi_{1}(x, t)=\phi_{1}(x)e^{-t}, \  \     \forall      \  x\in
R^{n}, \ t\geq 0.
\end{equation}
Then $ \Delta \psi_{1}= \psi_{1} $, and    $  \Box \psi_{1}=0 $.

One can see  \cite{Q. S. Zhang 1}, also see \cite{Q. S. Zhang 2}.

\section{ The proof of Theorem \ref{thm:1.1}} \vskip .5cm

The aim of this section is to prove Theorem 1.1, so we need to consider the following Cauchy problem
\begin{equation}\label{3.1}
\left \{
\begin{array}{lllll}
  \Box u(t,x ) =  u(t,x)u_{t} (t,x)^{ 2} + u^{ 4} (t,x),
  \   \   \
  x\in R^{ 2},  \   \    t>0,           \\
  t=0:  \     u=0,  \   \    \     \      u_{t}=\varepsilon g(|x|),   \  \
   \   x\in R^{ 2},
\end{array} \right.
\end{equation}
where
\begin{equation}\label{3.2}
   \    g(x)=g(|x|)> 0, \  \mbox{ for}  \   |x|<1      \  \   g(x)=0, \
 \   \mbox{ for}  \   |x|>1 \  \
\end{equation}
 We first prove that u is nonnegative.

 By the local existence of classical solutions, the solution to Cauchy problem \eqref{3.1} can be approximated by Picard iteration. Let
 $$u^{(0)}\equiv 0$$
 and
\begin{equation} \label{3.3}
\left \{
\begin{array}{lllll}
  \Box u^{(m)}(t,x ) =  u^{(m-1)}(t,x)u_{t}^{(m-1)} (t,x)^{ 2} + u^{(m-1)}(t,x)^{ 4} ,
  \   \   \
  x\in R^{ 2},  \   \    t>0,           \cr\noalign{\vskip2mm}
  t=0:  \     u^{(m)}=0,  \   \    \     \      u_{t}^{(m)}=\varepsilon g(|x|),   \  \
   \   x\in R^{ 2},
\end{array} \right.
\end{equation}
Then $\{u^{(m)}(t,x)\}$ is a series of approximate solutions to \eqref{3.1}.

Since $u^{(0)}\equiv 0$, by the positivity of the fundamental
solution of the wave operator in two space dimensions, we can prove
that all $u^{(m)}$ are nonnegative by induction. Let $m\rightarrow
\infty$, we can conclude that $ u$  is nonnegative.

The radial symmetric form of problem \eqref{3.1} can be written as
\begin{equation}\label{3.4}
\left \{
\begin{array}{lllll}
   u_{tt}-u_{rr}-\frac{u_r}{r} =  uu_{t}^{ 2} + u^{ 4},
  \   \  \   \    t>0,            \cr\noalign{\vskip2mm}
  t=0:  \     u=0,  \   \    \     \      u_{t}=\varepsilon g(r),   \  \

\end{array} \right.
\end{equation}
 where $$r=|x|,  \    \    \          x\in R^2.   $$
It follows from \eqref{3.4} that
\begin{equation} \label{3.5}
\left \{
\begin{array}{lllll}
   (\partial_{t}^2-\partial_{r}^2)\left(r^{\frac{1}{2}}u \right) = \frac{1}{4}r^{-\frac{3}{2}}u+ r^{\frac{1}{2}}uu_{t}^{ 2} +r^{\frac{1}{2}} u^{ 4},
  \   \  \   \    t>0,           \cr\noalign{\vskip2mm}
  t=0:  \     u=0,  \   \    \     \      r^{\frac{1}{2}}u_{t}=\varepsilon r^{\frac{1}{2}}g(r).   \  \

\end{array} \right.
\end{equation}
Let $G(r)=\frac{1}{2}r^{\frac{1}{2}}g(r)$, by D'Alembert's formula, in the domain $r>t$, we have
\begin{equation} \label{3.6}
\begin{array}{ll}
 r^{\frac{1}{2}}u(t,r) & =\varepsilon\displaystyle\int_{r-t}^{r+t}
  G(\lambda) d\lambda +\displaystyle\frac{1}{8}\displaystyle\int_{0}^{t} \int_{r-(t-\tau)}^{r+t-\tau}
 \frac{ u(\tau, \lambda)}{ \lambda^{ \frac{3}{2}} } d\lambda
 d\tau \cr\noalign{\vskip2mm}
  &  \   \   \   \      + \displaystyle\frac{1}{2}\displaystyle\int_{0}^{t} \int_{r-(t-\tau)}^{r+t-\tau}
 \lambda^{\frac{1}{2}}(uu_{t}^{ 2}+u^4)(\tau, \lambda) d\lambda d\tau.
 \end{array}
\end{equation}
Differentiate with respect to $ t $ yields:
\begin{equation} \label{3.7}
\begin{array}{lllll}
 r^{\frac{1}{2}}u_{t}(t,r) &=  \varepsilon G(t+r) + \varepsilon
 G(r-t)+ \displaystyle\frac{1}{8}\displaystyle\int_{0}^{t} \left[
 \frac{ u(\tau, \lambda)}{  \lambda^{ \frac{3}{2}} }
  \big |_{\lambda =r+t-\tau } + \frac{ u(\tau, \lambda)}{ \lambda^{ \frac{3}{2}} }
  \big |_{\lambda =r-(t-\tau) } \right]  d\tau \cr\noalign{\vskip3mm} & \   \
  \    \
+ \displaystyle\frac{1}{2}\displaystyle\int_{0}^{t} \left[
 \left( \lambda^{\frac{1}{2}}(uu_{t}^2+u^4)( \tau, \lambda )\right) \big|_{\lambda =r+t-\tau } +
   \left( \lambda^{\frac{1}{2}}(uu_{t}^2+u^4)( \tau, \lambda )  \right)  \big|_{\lambda =r-(t-\tau) } \right]d\tau.
 \end{array}
\end{equation}
Therefore, we obtain, in the domain $r>t$,
\begin{equation}\label{3.8}
r^{\frac{1}{2}}u(t,r)\ge \varepsilon\int_{r-t}^{r+t}G(\lambda
)d\lambda ,
\end{equation}
\begin{equation}\label{3.9}
r^{\frac{1}{2}}u_{t}(t,r)\ge \varepsilon G(t-r).
\end{equation}
It then easily follows that in the domain $t\ge \frac{1}{2}$ and $\frac{1}{4}\le r-t\le \frac{3}{4}$, we have
\begin{equation}\label{3.10}
u(t,r), u_t(t, r)\ge c_0\varepsilon r^{-\frac{1}{2}}.
\end{equation}
Let
\begin{equation}\label{3.11}
F(t)=\int_{R^2}u(t,x)dx,
\end{equation}
then by integrating \eqref{3.1}with respect to $ x $ , we obtain
$$F''(t)=\int_{R^2}(uu_t^2)(t,x)dx+\int_{R^2}(u^4)(t,x)dx .  $$
Thus, by the positivity of the solution $u$,  it follows
\begin{equation}\label{3.12}
F''(t)\ge \int_{R^2}(uu_t^2)(t,x)dx,
\end{equation}
and
\begin{equation}\label{3.13}
F''(t)\ge \int_{R^2}u^4(t,x)dx.
\end{equation}
Noting \eqref{3.10}, we obtain
$$
\begin{array}{ll}
 \displaystyle\int_{R^2}(uu_{t}^2)(t,x)dx &  =c\displaystyle\int_0^{t+1} (uu_{t}^2)(t,r)rdr
\cr\noalign{\vskip 3mm}   &  \ge
c\displaystyle\int_{t+\frac{1}{4}}^{t+\frac{3}{4}}
(uu_{t}^2)(t,r)rdr      \cr\noalign{\vskip 3mm}   &   \ge
c\varepsilon^3\displaystyle\int_{t+\frac{1}{4}}^{t+\frac{3}{4}}
r^{-\frac{1}{2}}dr
  \cr\noalign{\vskip 3mm}  &  \ge c\varepsilon^3(1+t)^{-\frac{1}{2}}.
 \end{array}$$
It then follows from \eqref{3.12}  that
\begin{equation}\label{3.15}
F(t)\ge c\varepsilon^3(1+t)^{\frac{3}{2}},\quad   t\ge 1.
\end{equation}
On the other hand, it follows from Holder's inequality that
$$F(t) =\int_{|x|\le t+1}u(t,x)dx
 \le \left(\int_{R^2}u^4(t,x)dx\right)^{\frac{1}{4}}\left(\int_{|x|\le t+1 }dx\right)^{\frac{3}{4}}
 \le C(1+t)^{\frac{3}{2}}\left(\int_{R^2}u^4(t,x)dx\right)^{\frac{1}{4}}.
$$
Noting \eqref{3.13}, we obtain
\begin{equation}\label {3.16}
F''(t)\ge C\frac{F^4(t)}{(1+t)^6}.
\end{equation}
Noting \eqref{3.15} and \eqref{3.16}, we may apply Lemma 2.1 (in which we take $\delta =\varepsilon ^3$, $\beta =4$, $a=\frac{3}{2}$, $\alpha =6$) to get the desired estimate of the life span
\begin{equation}\label{3.17}
T(\varepsilon )\le A\varepsilon^{-18},
\end{equation}
where $A $ is a positive constant which is independent of $
\varepsilon$.
    The proof of Theorem \ref{thm:1.1} is complete.

\section{ The proof of Theorem \ref{thm:1.2}} \vskip .5cm

Theorem \ref{thm:1.2} is  a consequence
 of  the blowup result and the upper bound  estimate about nonlinear
 differential inequalities in Lemma  \ref{thm:2.1}.

 We multiply the equation in  \eqref{1.3} by the test function $\psi_{1} (x, t)  \in
C^{2}(R^{n}\times R)$  and integrate over $R^{n} $,  then we use integration by parts.\\

 First,
$$\begin{array}{ll}
     \displaystyle \int_{R^{n }}
     \psi_{1}( u_{tt} -\Delta u )dx=  \displaystyle \int_{R^{n}}
     \psi_{1}|u_{ t}|^{p} dx +  \displaystyle \int_{R^{n}}
     \psi_{1}|u|^{q} dx .
\end{array}$$

 By the expression $\psi_{1}(x, t)=\phi_{1}(x)e^{-t} $,     we have
$$\begin{array}{ll}
     \displaystyle \int_{R^{n}}
     \psi_{1} \Delta u  dx=  \displaystyle \int_{R^{n}}
     \Delta  \psi_{1} u dx =   \displaystyle \int_{R^{n}}
     \psi_{1}  u dx .
\end{array}$$

So we have
\begin{equation}\label{4.1}
\begin{array}{ll}
     \displaystyle \int_{R^{n}}
      \left(   \psi_{1} u_{tt}  -   \psi_{1} u   \right)   dx=    \displaystyle \int_{R^{n}}
     \psi_{1}  ( |u_{ t}|^{p}+  |u |^{q}  ) dx .
\end{array}
\end{equation}

Notice that
\begin{equation}\label{4.2}
     \displaystyle\frac{d }{dt} \displaystyle\int_{R^{n}}
     \psi_{1}u_{t}dx= \displaystyle\int_{R^{n}} ( \psi_{1}\cdot u_{tt}-u_{t}\psi_{1})
     dx,
\end{equation}

\begin{equation} \label{4.3}
     \displaystyle\frac{d }{dt} \displaystyle\int_{R^{n }}
     (\psi_{1}u) dx= \displaystyle\int_{R^{n}}  \left[ (\psi_{1})_{t}\cdot u+u_{t}\cdot
     \psi_{1} \right]dx= \displaystyle\int_{R^{n}}
      \left[ \psi_{1}\cdot u_{t}-u \psi_{1}\right]dx.
\end{equation}

Adding up the above two expressions, we obtain the following
\begin{equation} \label{4.4}
     \displaystyle\frac{d }{dt} \displaystyle\int_{R^{n}}
     \left(\psi_{1}u_{t}+ \psi_{1}u \right) dx= \displaystyle\int_{R^{n}}
  \left( \psi_{1}\cdot u_{tt}-u\cdot
     \psi_{1} \right) dx=  \displaystyle\int_{R^{n}} \left (  \psi_{1} \cdot |u_{t}|^{p}
       + \psi_{1} \cdot | u |^{ q}  \right) dx.
\end{equation}

So we have
\begin{equation} \label{4.5}
\begin{array}{ll}
    &   \displaystyle\int_{R^{n}}
     \left(\psi_{1}u_{t}+ \psi_{1}u \right) dx \cr\noalign{\vskip 3mm}  &= \displaystyle\int_{R^{n}}
      \left(\psi_{1}u_{t}+ \psi_{1}u \right) dx|_{t=0}+
\displaystyle\int_{0}^{t}\int_{ R^{n}}
   \psi_{1}\cdot    \left (  |u_{t}|^{p}
       +  | u | ^{ q}  \right)
    dxd \tau \cr\noalign{\vskip 3mm} &
=   \varepsilon \displaystyle\int_{R^{n}}
      \phi_{1}(x)\left[f(x)+ g(x) \right] dx+ \displaystyle\int_{0}^{t}\int_{R^{n}}
   \psi_{1}\cdot   \left (  |u_{t}|^{p}
       +  | u |^{ q}  \right)  dxd \tau.
 \end{array}
\end{equation}

Adding the expressions   \eqref{4.1}  and   \eqref{4.5}, we have
\begin{equation} \label{4.6}
\begin{array}{ll}
     &  \displaystyle\int_{R^{n}}
     \left(\psi_{1}u_{t}+ \psi_{1}u_{ tt} \right) dx
  \cr\noalign{\vskip 3mm}   & =
  \varepsilon \displaystyle\int_{R^{n}}
      \phi_{1}(x)\left[f(x)+ g(x) \right] dx  +  \displaystyle \int_{R^{n}}
     \psi_{1}  ( |u_{ t}|^{p}+  |u | ^{q}  ) dx + \displaystyle\int_{0}^{t}\int_{R^{n}}
   \psi_{1}\cdot   \left (  |u_{t}|^{p}
       +  |u |^{ q}  \right)  dxd \tau
     \cr\noalign{\vskip 3mm}   &
 \geq    \varepsilon \displaystyle\int_{R^{n}}
      \phi_{1}(x)\left[f(x)+ g(x) \right] dx    =  C_{0} \varepsilon.
 \end{array}
\end{equation}

 Also, we know that
\begin{equation} \label{4.7}
\begin{array}{ll}
     \displaystyle\frac{d }{dt} \displaystyle\int_{R^{n}}
     \psi_{1}u_{t} dx + 2 \displaystyle\int_{R^{n}}
     \psi_{1}\cdot u_{t}dx&= \displaystyle\int_{R^{n}}
  \left[ \psi_{1}u_{tt}+u_{t}(\psi_{1})_{t}+2\psi_{1}u_{t} \right]
  dx \cr\noalign{\vskip 3mm} &= \displaystyle\int_{R^{n}} \left( \psi_{1}u_{tt}+ \psi_{1}u_{t}
  \right)dx   \cr\noalign{\vskip 3mm} & \geq  C_{0} \varepsilon .
 \end{array}
\end{equation}

  Multiplying the above differential inequality by $e^{2t}$,  we get the following expression
    \begin{equation} \label{4.8}
     \displaystyle\frac{d }{dt} \left( e^{2t}   \displaystyle\int_{R^{n}}
     \psi_{1}u_{t} dx  \right)\geq  C_{0} \varepsilon e^{2t} .
      \end{equation}

So we have
  \begin{equation} \label{4.9}
    e^{2t}   \displaystyle\int_{R^{n}}
     \psi_{1}u_{t} dx \geq      C_{0}  ( e^{2t} -1)     \varepsilon + \varepsilon  \displaystyle\int_{R^{n}} \phi_{ 1} gdx  .
 \end{equation}

Therefore,   noting the positivity of $\phi_{ 1}$ and $g$,    we
have
 \begin{equation} \label{4.10}       \displaystyle\int_{R^{n}}
     \psi_{1}u_{t} dx \geq      \widetilde{C_{0}}  \varepsilon ,
     \end{equation}
 where $\widetilde{C_{0}} $ is a positive constant.

   Let  $  F(t) =      \displaystyle\int_{R^{n}}  u dx,
    $   we integrate the equation   $$    \Box u =  |  u_{t}|^{ p} +  | u |^{ q}  ,   $$
 we have
  \begin{equation} \label{4.11}
   \partial_{ t}^{ 2}  \displaystyle\int_{R^{n}}  u dx \geq
\displaystyle\int_{R^{n}}|  u_{t}|^{ p} dx +
              \displaystyle\int_{R^{n}}  | u |^{ q} dx .    \end{equation}
That is
   \begin{equation} \label{4.12}
     F''(t)  \geq    \displaystyle\int_{R^{n}}  |  u_{t}|^{ p} dx  +   \displaystyle\int_{R^{n}}  |u |^{ q}
     dx.
      \end{equation}

 By Holder's inequality, we can obtain
    \begin{equation} \label{4.13}       \displaystyle\int_{R^{n}}    u_{t}\psi_{1} dx \leq
   \left(  \displaystyle\int_{ | x|\leq t+1   }  |  u_{t}|^{ p} dx   \right)^{  \frac{1}{ p}}
   \left(  \displaystyle\int_{ | x|\leq t+1   }  |  \psi_{1} |^{ p'} dx   \right)^{  \frac{1}{ p' }}  ,
     \end{equation}
where $ p $ and $ p'$ satisfies $ \frac{ 1}{ p} + \frac{ 1}{ p'}=1
$.

 Noting \eqref{2.2}, it turns out that
     \begin{equation} \label{4.14}
     \begin{array}{ll}
   \left(  \displaystyle\int_{ | x|\leq t+1   }  |  \psi_{1} |^{ p'} dx   \right)^{  \frac{1}{ p'}}
   &  =  e^{ -t}  \left(  \displaystyle\int_{ | x|\leq t+1   }  |  \phi_{1}(x) |^{ p'} dx   \right)^{  \frac{1}{ p'}}
   \cr\noalign{\vskip 3mm} &
   \leq  C e^{ -t} \left(  \displaystyle\int_{  0}^{ t+1}    e^{  p'r}
   (1+r)^{  -\frac{ ( n-1 ) p'}{2 }} r^{ n-1 } dr   \right)^{  \frac{1}{ p'}}
  \cr\noalign{\vskip 3mm} &
    \leq  C e^{ -t} \left[ ( t+1 )^{  n-1- \frac{ ( n-1 ) p'}{ 2}}    \right]^{  \frac{ 1}{ p' }}
     \left(    \displaystyle\int_{  0}^{ t+1}    e^{  p' r}
    dr   \right)^{  \frac{1}{ p'}}
  \cr\noalign{\vskip 3mm} &
    =  C e^{ -t}  ( t+1 )^{   \frac{n- 1}{ p'} -\frac{ n-1}{ 2}  }
     \left[   \frac{ 1}{p' } e^{  p' (t+1) }  -  \frac{ 1}{ p' }     \right]^{  \frac{1}{ p'}}
 \cr\noalign{\vskip 3mm} &
  \leq   C e^{ -t}  ( t+1 )^{    \frac{n- 1}{ p'} -\frac{ n-1}{ 2} }
        ( \frac{ 1}{ p' } )^{\frac{ 1 }{ p' } }   e^{  t+1 }
     \cr\noalign{\vskip 3mm} &    =  C(1+t)^{   \frac{n- 1}{ p'} -\frac{ n-1}{ 2}   }
   .
    \end{array}
    \end{equation}

Noting  \eqref{4.10} and  \eqref{4.14},   it follows from
\eqref{4.13} that
  \begin{equation} \label{4.15}
         \displaystyle\int_{R^{n}}  |  u_{t}|^{ p} dx   \geq
  C   \displaystyle \frac{    (   \widetilde{C_{0}}  \varepsilon    )^{ p}  }{ (1+t)^{  ( \frac{n- 1}{ p'} -\frac{ n-1}{ 2}  ) p    }  }
  .
    \end{equation}
Since $  p'= \frac{ p}{ p-1} $, so we have
$$  \left [  \frac{n- 1}{ p'} -\frac{ n-1}{ 2}   \right]  p = \frac{ (n- 1) p }{ p'} -\frac{ (n-1) p}{ 2} =
\frac{ (n-1) (p-2)}{ 2},            $$
  thus the expression
\eqref{4.15} leads to the following
 \begin{equation} \label{4.16}
         \displaystyle\int_{R^{n}}  |  u_{t}|^{ p} dx   \geq
  C   \displaystyle \frac{    (   \widetilde{C_{0}}  \varepsilon    )^{ p}  }{ (1+t)^{ \frac{ (n-1) (p-2)}{ 2}     }  }
  .
    \end{equation}

By Holder's inequality, we obtain
  \begin{equation} \label{4.17}
    \begin{array}{ll}
     F(t) =    \displaystyle\int_{R^{n}}  u dx   & \leq
   \left(  \displaystyle\int_{R^{n}}   | u |^{ q } dx   \right)^{  \frac{1}{ q}}
    \left(  \displaystyle\int_{ | x|\leq  t+1  }     dx   \right)^{  \frac{1}{ q'}}
   \cr\noalign{\vskip 3mm} &
      \leq  C  \left(  \displaystyle\int_{R^{n}}    |u |^{ q} dx   \right)^{  \frac{1}{ q}}
    (1+t)^{ \frac{ n}{  q'}},
 \end{array}
     \end{equation}
 where  $ q $ and $ q'$  are conjugate numbers, they satisfy $ \frac{ 1}{ q} + \frac{ 1}{ q'}=1
$.
 Therefore, we have
 \begin{equation} \label{4.18}       \displaystyle\int_{R^{n}}
 |u |^{ q } dx  \geq  C \displaystyle \frac{  F(t)^{ q} }{ (1+t)^{ \frac{ nq}{ q'}} } .
         \end{equation}
Noticing that  $ q'= \frac{ q }{ q-1} $, so  $ \frac{ nq}{ q'} =
n(q-1) $.
 So the expression   \eqref{4.18} leads to the following
\begin{equation} \label{4.19}       \displaystyle\int_{R^{n}}
 |u |^{ q } dx  \geq  C \displaystyle \frac{  F(t)^{ q} }{ (1+t)^{ n( q-1 )} } .
         \end{equation}
  Hence, $ F(t) $  satisfies the following inequality
 \begin{equation} \label{4.20}  F''(t)    \geq  C \displaystyle
\frac{  F(t)^{ q} }{ (1+t)^{ n(q-1 )} } .
     \end{equation}

 On the other hand,  by \eqref{4.16}, we get:
  \begin{equation} \label{4.21}
   F''(t)    \geq  C \displaystyle \frac{  \varepsilon^{ p} }{ (1+t)^{
      \frac{ (n-1) (p-2)}{ 2}       } } ,
     \end{equation}

 So  integrating the above expression twice,  we have the following
    \begin{equation} \label{4.22}    F(t) \geq C \varepsilon^{ p}(1+t)^{  2-\frac{ (n-1) (p-2)}{ 2}     }
    .
     \end{equation}

 We take $  a=    2-\frac{ (n-1) (p-2)}{ 2}   $,  $ \alpha = n(q-1)  $, $ \beta=q>1 $ in Lemma  \ref{thm:2.1},  from
the conditions \eqref{1.7} and  \eqref{1.8} in theorem
\ref{thm:1.2}, we have   $ \beta>1$, $ a\geq 1$  can be deduced from \eqref{1.7} and  $ ( \beta -1 ) a>\alpha -2 $  can be deduced from
\eqref{1.8}.  So by the Lemma  \ref{thm:2.1},  $F(t)$ will blow up
in finite time and thus  the solutions to problem   \eqref{1.3} will
blow up in finite time,  and also  we have the following
   \begin{equation} \label{4.23}
    \begin{array}{ll}
   T(\delta) &  \leq c \delta^{ -\frac{(\beta-1)}{ (\beta-1)  a-\alpha +2}
   }     \cr\noalign{\vskip 3mm} &
    =  c \delta^{ - \frac{ (q-1)   }{ q+1-\frac{ (n-1)p(q-1)}{2}   } },

    \end{array}
    \end{equation}
here we take $ \delta = \varepsilon^{ p} $, then we have the
estimate for the lifespan of the solution to the problem
\eqref{1.3}:
  \begin{equation} \label{4.24}   T( \varepsilon) \leq
   C\varepsilon^{ - \frac{ p(q-1)   }{q+1-\frac{ (n-1)p(q-1)}{2}  } }
    = C\varepsilon^{-  \frac{ 2p(q-1)   }{ 2q+2- (n-1) p(q-1)  }  },    \end{equation}
 where $C $ is a positive constant which is independent of $ \varepsilon$.
    The proof of Theorem \ref{thm:1.2} is complete.

\section*{Acknowledgments.}
The  authors would like to thank Professor Takamura for  pointing
out the reference Katayama  \cite{Katayama} to us and for helpful
discussions.

   Yi Zhou's
research is supported by  the Key Laboratory of Mathematics for
Nonlinear Sciences,  Ministry of Education of China and Shanghai Key
Laboratory for Contemporary and Applied Mathematics, and NSFC(
No.11031001 and 11121101), and 111 project;
 Wei Han's research is  supported by the Youth Science
Foundation of Shanxi Province (2010021001-2), the National Sciences
Foundation of China(10901145),  the Top Young Academic Leaders of
Higher Learning Institutions of Shanxi.

\end{document}